\documentclass[a4paper]{amsart}

\usepackage{amssymb}
\usepackage{color}
\usepackage{epsf}
\usepackage{graphicx}
\usepackage{bbm}

\newtheorem{theorem}{Theorem}[section]
\newtheorem{lem}[theorem]{Lemma}
\newtheorem{prop}[theorem]{Proposition}

\newtheorem{coro}[theorem]{Corollary}
\newtheorem{fact}[theorem]{Fact}

\theoremstyle{definition}
\newtheorem{defi}[theorem]{Definition}
\newtheorem{ex}[theorem]{Example}
\newtheorem{rem}[theorem]{Remark}

\numberwithin{theorem}{section}

\newcommand{\OO}{\mathcal{O}}
\newcommand{\oo}{\thinspace\scriptstyle{\mathcal{O}}\displaystyle}
\newcommand{\zetan}{\zeta_{n}}
\newcommand{\zetam}{\zeta_{m}}
\newcommand{\Kn}{\mathbbm{K}_{n}}
\newcommand{\kn}{\mathbbm{k}_{n}}
\newcommand{\Zn}{\mathbbm{Z}[\zeta_{n}]}

\newcommand{\Qn}{\mathbbm{Q}(\zeta_{n})}

\newcommand{\LN}{\mathbbm{N}}
\newcommand{\LC}{\mathbbm{C}}
\newcommand{\C}{\mathbbm{C}}
\newcommand{\LQ}{\mathbbm{Q}}
\newcommand{\LR}{\mathbbm{R}}
\newcommand{\LZ}{\mathbbm{Z}}
\newcommand{\R}{\mathbbm{R}}
\newcommand{\Q}{\mathbbm{Q}}
\newcommand{\Z}{\mathbbm{Z}}

\newdimen{\standardlabelwidth}
\standardlabelwidth-\standardlabelwidth
  \advance\standardlabelwidth by 2\normalparindent
\newcommand{\standardlabel}[1]{#1\kern\standardlabelwidth}

\begin{document}


\title[A note on affinely regular polygons]{A note on affinely regular polygons}
\author{Christian Huck}
\address{Department of Mathematics and Statistics\\
  The Open University\\ Walton Hall\\ Milton Keynes\\ MK7 6AA\\
   United Kingdom}
\email{c.huck@open.ac.uk}
\thanks{The author was supported by the German Research Council (Deutsche Forschungsgemeinschaft), within the CRC 701, and by EPSRC via Grant EP/D058465/1.}

\begin{abstract}
The affinely regular polygons in certain planar sets are characterized. It is also shown that the obtained results apply to cyclotomic model sets and, additionally, have consequences in the discrete tomography of these sets.    
\end{abstract}

\maketitle


\section{Introduction}
Chrestenson~\cite{C} has shown that any
(planar) regular polygon whose vertices are contained in
$\mathbbm{Z}^{d}$ for some $d\geq 2$ must have $3,4$ or $6$
vertices. More generally, Gardner and Gritzmann~\cite{GG} have
characterized the numbers of vertices of affinely regular \emph{lattice polygons}, i.e.,
images of non-degenerate regular polygons under a non-singular affine
transformation of the plane whose vertices are
contained in the square lattice $\LZ^2$ or, equivalently, in some arbitrary planar lattice $L$. It turned out that the affinely regular lattice polygons are precisely
the affinely regular triangles, parallelograms and hexagons. As a first step beyond the case of planar lattices, this short text provides a generalization of this result to planar sets $\varLambda$ that  are non-degenerate in some sense and satisfy a certain affinity condition on finite scales (Theorem~\ref{charac}). The obtained characterization can be expressed in terms of a simple inclusion of real field extensions of $\LQ$ and particularly applies to \emph{algebraic Delone sets}, thus including \emph{cyclotomic model sets}. These cyclotomic model sets range from periodic examples, given by the vertex sets of the square tiling and the triangular tiling, to aperiodic examples like the vertices of the Ammann-Beenker tiling, of the T\"ubingen triangle tiling and of the shield tiling, respectively. I turns out that, for cyclotomic model sets $\varLambda$, the numbers of vertices of affinely regular polygons in $\varLambda$ can be characterized by a simple divisibility condition (Corollary~\ref{th1mod}). In particular, the result on affinely regular lattice polygons is contained as a special case (Corollary~\ref{cor2}(a)). Additionally, it is shown that the obtained divisibility condition implies a weak estimate in the \emph{discrete tomography} of cyclotomic model sets (Corollary~\ref{corofin}).

\section{Preliminaries and notation}\label{sec1}
 
Natural numbers are always assumed to be positive, i.e., 
$\mathbbm{N}\,=\,\{1,2,3,\dots\}$ and we denote by $\mathcal{P}$ the set of rational primes. If
$k,l\in \mathbbm{N}$, then $\operatorname{gcd}(k,l)$ and $\operatorname{lcm}(k,l)$
 denote their greatest common divisor and least common multiple, respectively. The group of units of a given ring $R$ is denoted by $R^{\times}$. As usual, for a complex number $z\in\LC$, $\vert z\vert$ denotes the complex absolute value, i.e., $\vert z\vert=\sqrt{z\bar{z}}$, where $\bar{.}$ denotes the complex
conjugation. The unit circle in
$\mathbbm{R}^{2}$ is denoted by $\mathbb{S}^{1}$, i.e.,
 $\mathbb{S}^{1}=\{x\in\R^2\,|\,\vert x \vert =1\}$. Moreover, the elements of $\mathbb{S}^{1}$ are also called
{\em directions}. For $r>0$ and $x\in\R^{2}$,
$B_{r}(x)$ denotes the open ball of radius $r$ about $x$. A subset $\varLambda$ of the plane is called {\em uniformly discrete} if there is a radius
$r>0$ such that every ball $B_{r}(x)$ with $x\in\mathbbm{R}^{2}$ contains at most one point of
  $\varLambda$. Further, $\varLambda$ is called {\em relatively dense} if there is a radius $R>0$
  such that every ball $B_{R}(x)$ with 
  $x\in\mathbbm{R}^{2}$ contains at least one point of $\varLambda$. $\varLambda$ is called a {\em Delone set} (or {\em Delaunay set}) if it is both uniformly
  discrete and relatively dense. For a subset $S$ of the plane, we denote by
$\operatorname{card}(S)$, $\mathcal{F}(S)$, $\operatorname{conv}(S)$ and $\mathbbm{1}_{S}$ the cardinality, set of finite subsets, convex hull and
characteristic function of $S$, respectively. A direction
$u\in\mathbb{S}^{1}$ is called an $S${\em-direction} if it is
parallel to a non-zero element of the difference set
$S-S:=\{s-s'\,|\,s,s'\in S\}$ of $S$. Further, a finite subset $C$ of $S$ is called a {\em convex subset of} $S$ if its
convex hull contains no new points of $S$, i.e., if $C =
\operatorname{conv}(C)\cap S$ holds. Moreover, the set of all convex subsets
of 
$S$ is denoted by $\mathcal{C}(S)$. Recall that a {\em linear transformation} (resp., {\em affine transformation})
$\Psi\!:\, \mathbbm{R}^{2} \rightarrow \mathbbm{R}^{2}$ of the Euclidean plane is given
by $z \mapsto Az$ (resp., $z \mapsto Az+t$), where $A$ is a
real $2\times 2$ matrix and $t\in \mathbbm{R}^{2}$. In both cases,
$\Psi$ is called {\em singular} when $\operatorname{det}(A)= 0$;
otherwise, it is non-singular. A {\em homothety} $h\!:\, \mathbbm{R}^{2} \rightarrow
\mathbbm{R}^{2}$ is given by $z \mapsto \lambda z + t$, where
$\lambda \in \R$ is positive and $t\in \mathbbm{R}^{2}$. A {\em convex polygon} is the convex hull of a finite set of points in $\R^2$. For a subset $S \subset \R^2$, a {\em polygon in} $S$ is a convex polygon with all vertices in $S$. A {\em regular polygon} is always assumed to be planar, non-degenerate and convex. An {\em affinely regular polygon} is a non-singular affine image of a regular polygon. In particular, it must have at least $3$ vertices. Let $U\subset \mathbb{S}^{1}$ be a finite set of directions. A non-degenerate convex polygon $P$ is called a $U${\em -polygon} if it has the property that whenever $v$ is a vertex of $P$ and $u\in U$, the line $\ell_{u}^{v}$ in the plane in direction $u$ which passes through $v$ also meets another vertex $v'$ of $P$. For a subset $\varLambda\subset\LC$, we denote by $\mathbbm{K}_{\varLambda}$ the intermediate field of $\LC/\LQ$ that is given by 
$$
\mathbbm{K}_{\varLambda}\,\,:=\,\,\Q\left(\big(\varLambda-\varLambda\big)\cup\big(\overline{\varLambda-\varLambda}\big)\right)\,,
$$
where $\varLambda-\varLambda$ denotes the difference set of $\varLambda$. Further, we set 
$
\mathbbm{k}_{\varLambda}:=\mathbbm{K}_{\varLambda}\cap\R$, the maximal real subfield of $\mathbbm{K}_{\varLambda}$.

\begin{rem}
Note that $U$-polygons have an even number of vertices. Moreover, an affinely regular polygon with an even number of vertices is a $U$-polygon if and only if each direction of $U$ is parallel to one of its edges.     
\end{rem}

For $n\in
\mathbbm{N}$, we always let $\zetan := e^{2\pi i/n}$, as a specific
choice for a
primitive $n$th root of unity in $\LC$. Let $\Q(\zetan)$ be
the corresponding cyclotomic field. It is well known that
$\Q(\zetan+\bar{\zeta}_{n})$ is the maximal real subfield of
$\Q(\zetan)$; see~\cite{Wa}. Throughout this text, we shall use the
notation $$\mathbbm{K}_{n}=\Qn,\;
\mathbbm{k}_{n}=\Q(\zetan+\bar{\zeta}_{n}),\; \OO_{n}=\Zn,\; \oo_{n}
=\Z[\zetan+\bar{\zeta}_{n}]\,.$$ Except for the one-dimensional cases $\mathbbm{K}_{1}=\mathbbm{K}_{2}=\LQ$, $\mathbbm{K}_{n}$ is an imaginary extension of $\LQ$. Further, $\phi$ will always denote Euler's phi-function, i.e., $$\phi(n) =
\operatorname{card}\left(\big\{k \in \mathbbm{N}\, |\,1 \leq k \leq n
  \textnormal{ and } \operatorname{gcd}(k,n)=1\big\}\right)\,.$$ Occasionally, we
identify $\C$ with $\R^{2}$. Primes $p\in\mathcal{P}$ for which the number $2p+1$ is prime as well are called Sophie Germain prime numbers. We denote by $\mathcal{P}_{\rm SG}$ the set of Sophie Germain prime numbers.  They are the primes $p$ such that the equation
$\phi(n)=2p$
has solutions. It is not known whether there are infinitely many
Sophie Germain primes. The first few are
\begin{eqnarray*}&&\{2,3,5,11,23,29,41,53,83,89,113,131,173,\\
  &&\hphantom{\{}179,191,233,239,251,
           281,293,359,419,\dots\}\,,\end{eqnarray*} 
see entry A005384 of~\cite{Sl} for further details. We need the following facts from the theory of cyclotomic fields.

\begin{fact}[Gau\ss]\cite[Theorem 2.5]{Wa}\label{gau}
$[\mathbbm{K}_{n} : \mathbbm{Q}] = \phi(n)$. The field extension $\mathbbm{K}_{n}/ \mathbbm{Q}$
is a Galois extension with Abelian Galois group $G(\mathbbm{K}_{n}/
\mathbbm{Q}) \simeq (\Z / n\Z)^{\times}$,
where $a\, (\textnormal{mod}\, n)$ corresponds to the automorphism given by\/ $\zetan \mapsto \zetan^{a}$. 
\end{fact}

Since $\mathbbm{k}_n$ is the maximal real subfield of the $n$th cyclotomic field $\mathbbm{K}_{n}$, Fact~\ref{gau} immediately gives the following result. 

\begin{coro}\label{cr5}
If $n\geq 3$, one has $[\mathbbm{K}_{n} : \mathbbm{k}_{n}]= 2$. Thus, a $\mathbbm{k}_{n}$-basis of $\mathbbm{K}_{n}$ is given by $\{1,\zeta_n\}$. The field extension $\mathbbm{k}_{n} / \mathbbm{Q}$  is a Galois extension with Abelian
Galois group $G(\mathbbm{k}_{n} / \mathbbm{Q})\simeq (\Z /
n\Z)^{\times}/\{\pm 1\, (\textnormal{mod}\, n)\}$ of order $[\mathbbm{k}_{n} : \mathbbm{Q}] =
\phi(n)/2$.
\end{coro}

Consider an
algebraic number field $\mathbbm{K}$, i.e., a finite extension of $\Q$. A full $\mathbbm{Z}$-module $\mathcal{O}$ in $\mathbbm{K}$ (i.e., a free $\mathbbm{Z}$-module of rank $[\mathbbm{K}:\LQ]$) which contains the number $1$ and
is a ring is called an \emph{order} of $\mathbbm{K}$. Note that every $\LZ$-basis of $\mathcal{O}$ is simultaneously a $\LQ$-basis of $\mathbbm{K}$, whence $\Q\mathcal{O}=\mathbbm{K}$ in particular. It turns out
that among the various orders of $\mathbbm{K}$ there is one
\emph{maximal order} which contains all the other orders, namely the
ring of integers $\mathcal{O}_{\mathbbm{K}}$ in $\mathbbm{K}$; see~\cite[Chapter 2, Section 2]{Bo}. For cyclotomic fields, one has the
following well-known result.

\begin{fact}\cite[Theorem 2.6 and Proposition 2.16]{Wa}\label{p1}
For $n\in \mathbbm{N}$, one has:
\begin{itemize}
\item[(a)]
$\OO_{n}$ is the ring of cyclotomic integers in $\mathbbm{K}_{n}$, and hence its maximal order.
\item[(b)]
$\oo_{n}$ is the ring of integers in $\mathbbm{k}_{n}$, and hence its maximal order.
\end{itemize}
\end{fact}

\begin{lem}\label{cyclosec}
If $m,n \in \mathbbm{N}$, then $\mathbbm{K}_{m} \cap \mathbbm{K}_{n} = \mathbbm{K}_{\operatorname{gcd}(m,n)}$.
\end{lem}
\begin{proof}
The assertion follows from similar arguments as in the proof of the special case $(m,n)=1$; compare \cite[Ch. VI.3, Corollary 3.2]{La}. Here, one has to observe $\mathbbm{Q}(\zeta_{m},\zetan)=\mathbbm{K}_{m}\mathbbm{K}_{n}= \mathbbm{K}_{\operatorname{lcm}(m,n)}$ and then to employ the identity \begin{equation}\label{eq2}\phi(m)\phi(n)=\phi(\operatorname{lcm}(m,n))\phi(\operatorname{gcd}(m,n))\end{equation} instead of merely using the multiplicativity of the arithmetic function $\phi$. 
\end{proof}

\begin{lem}\label{incl}
Let $m,n \in \mathbbm{N}$. The following statements are equivalent:
\begin{itemize}
\item[(i)]
$\mathbbm{K}_{m} \subset \mathbbm{K}_{n}$.
\item[(ii)]
$m|n$, or $m\; \equiv \;2  \;(\operatorname{mod} 4)$ and $m|2n$.
\end{itemize}
\end{lem}
\begin{proof}
For direction (ii) $\Rightarrow$ (i), the assertion is clear if 
$m|n$. Further, if $m\; \equiv \;2
\;(\operatorname{mod} 4)$, say $m=2o$ for a suitable odd number $o$,
and $m|2n$, then $\mathbbm{K}_{o} \subset \mathbbm{K}_{n}$ (due to
$o|n$). However, Fact~\ref{gau} shows that the inclusion of fields
$\mathbbm{K}_{o} \subset \mathbbm{K}_{2o}=\mathbbm{K}_{m}$ cannot be
proper since we have, by means of the multiplicativity of $\phi$, the equation 
$\phi(m)=\phi(2o)=\phi(o)$. This gives $\mathbbm{K}_{m}\subset \mathbbm{K}_{n}$. 

For direction (i) $\Rightarrow$ (ii), suppose $\mathbbm{K}_{m} \subset \mathbbm{K}_{n}$. Then, Lemma~\ref{cyclosec} implies $\mathbbm{K}_{m} = \mathbbm{K}_{\operatorname{gcd}(m,n)}$, whence \begin{equation}\label{equation}\phi(m)=\phi(\operatorname{gcd}(m,n))\,\end{equation} by Fact~\ref{gau} again. Using the multiplicativity of $\phi$ together with $\phi(p^{j})\,=\,p^{j-1}\,(p-1)$ for $p\in\mathcal{P}$ and $j\in\mathbbm{N}$, we see that, given the case $\operatorname{gcd}(m,n)< m$, Equation (\ref{equation}) can only be fulfilled if $m\; \equiv \;2  \;(\operatorname{mod} 4)$ and $m|2n$. The remaining case $\operatorname{gcd}(m,n)= m$ is equivalent to the relation $m|n$.
\end{proof}

\begin{coro}\label{unique}
Let $m,n \in \mathbbm{N}$. The following statements are equivalent:
\begin{itemize}
\item[(i)]
$\mathbbm{K}_{m} = \mathbbm{K}_{n}$.
\item[(ii)]
$m=n$, or $m$ is odd and $n=2m$, or $n$ is odd and $m=2n$.
\end{itemize}
\end{coro}

\begin{rem}\label{unrem}
Corollary~\ref{unique} implies that, for $m,n\; \not\equiv \;2
\;(\operatorname{mod} 4)$, one has the identity $\mathbbm{K}_{m} = \mathbbm{K}_{n}$ if and only if $m=n$.
\end{rem}

\begin{lem}\label{unique2}
Let $m,n \in \mathbbm{N}$ with $m,n \geq 3$. Then, one has:
\begin{itemize}
\item[(a)]
$\mathbbm{k}_{m}=\mathbbm{k}_{n} \;\Leftrightarrow\; \mathbbm{K}_{m}=\mathbbm{K}_{n} \mbox{ or\, } m,n \in \{3,4,6\}.$
\item[(b)]
$\mathbbm{k}_{m}\subset\mathbbm{k}_{n} \;\Leftrightarrow\; \mathbbm{K}_{m} \subset \mathbbm{K}_{n} \mbox{ or\, } m\in \{3,4,6\}.$
\end{itemize}
\end{lem}
\begin{proof}
For claim (a), let us suppose $\mathbbm{k}_{m}=\mathbbm{k}_{n}=:\mathbbm{k}$ first. Then, Fact~\ref{gau} and Corollary~\ref{cr5} imply that $[\mathbbm{K}_{m}:\mathbbm{k}]=[\mathbbm{K}_{n}:\mathbbm{k}]=2$. Note that $\mathbbm{K}_{m} \cap \mathbbm{K}_{n} = \mathbbm{K}_{\operatorname{gcd}(m,n)}$ is a cyclotomic field containing $\mathbbm{k}$. It follows that either $\mathbbm{K}_{m} \cap \mathbbm{K}_{n}=\mathbbm{K}_{\operatorname{gcd}(m,n)}=\mathbbm{K}_{m}=\mathbbm{K}_{n}$ or $\mathbbm{K}_{m} \cap \mathbbm{K}_{n}=\mathbbm{K}_{\operatorname{gcd}(m,n)}= \mathbbm{k}$ and hence $\mathbbm{k}_{m}=\mathbbm{k}_{n}= \mathbbm{k}= \mathbbm{Q}$, since the latter is the only real cyclotomic field. Now, this implies $m,n\in\{3,4,6\}$; see also Lemma~\ref{phin2p}(a) below. The other direction is obvious. Claim (b) follows immediately from the part (a).   
\end{proof}

\begin{lem}\label{phin2p}
Consider $\phi$ on $\{n\in \mathbbm{N}\, |\, n\; \not\equiv \;2  \;(\operatorname{mod} 4)\}$. Then, one has: 
\begin{itemize}
\item[(a)] $\phi(n)/2=1$ if and only if $n\in\{3,4\}$.
\item[(b)] $\phi(n)/2 \in \mathcal{P}$ if and only if $n \in \mathcal{S}:=\{8,9,12\} \cup \{2p+1\, | \, p \in \mathcal{P}_{\rm SG}\}$.
\end{itemize}
\end{lem}
\begin{proof}
The equivalences follow from the multiplicativity of $\phi$ in conjunction with the identity $\phi(p^{j})\,=\,p^{j-1}\,(p-1)$ for $p\in\mathcal{P}$ and $j\in\mathbbm{N}$.
\end{proof}

\begin{rem}
Let $n\; \not\equiv \;2  \;(\operatorname{mod} 4)$. By
Corollary~\ref{cr5}, for $n\geq 3$, the field extension $\mathbbm{k}_{n} / \mathbbm{Q}$ is a Galois
extension with Abelian Galois group $G(\mathbbm{k}_{n} / \mathbbm{Q})$
of order $\phi(n)/2$. Using Lemma~\ref{phin2p}, one sees that
$G(\mathbbm{k}_{n} / \mathbbm{Q})$ is trivial if and only if
$n\in\{1,3,4\}$, and simple if and only if $n\in\mathcal{S}$, with
$\mathcal{S}$ as defined in Lemma~\ref{phin2p}(b). 
\end{rem}

\section{The characterization}\label{sec2}

The following notions will be of crucial importance. 

\begin{defi}\label{algdeldef}
For a set $\varLambda\subset\R^2$, we define the following properties:
\begin{itemize}
\item[(Alg)]
$\left[\mathbbm{K}_{\varLambda}:\Q\right]<\infty$.
\item[(Aff)]
For all $F\in\mathcal{F}(\mathbbm{K}_{\varLambda})$, there is a
non-singular affine transformation $\Psi\!:\, \mathbbm{R}^{2} \rightarrow
\mathbbm{R}^{2}$ such that $h(F)\subset \varLambda$.
\end{itemize} 
Moreover, $\varLambda$ is called \emph{degenerate} when $\mathbbm{K}_{\varLambda}\subset \LR$; otherwise, $\varLambda$ is non-degenerate.
\end{defi} 

\begin{rem}\label{mrs}
If $\varLambda\subset\R^2$ satisfies property~(Alg), then one has $\left[\mathbbm{k}_{\varLambda}:\Q\right]<\infty$, i.e., $\mathbbm{k}_{\varLambda}$ is a real algebraic number field.
\end{rem}

Before we turn to examples of planar sets $\varLambda$ having properties (Alg) and (Aff), let us prove the central result of this text, where we use arguments similar to the ones used by Gardner and Gritzmann in the proof of~\cite[Theorem 4.1]{GG}.

\begin{theorem}\label{charac}
Let $\varLambda\subset\R^2$ be non-degenerate with property~(Aff). Further, let $m\in\mathbbm{N}$ with $m\geq 3$. The following statements are equivalent:
\begin{itemize}
\item[(i)]
There is an affinely regular $m$-gon in $\varLambda$.
\item[(ii)]
$\mathbbm{k}_{m}\subset\mathbbm{k}_{\varLambda}$.
\end{itemize}
If $\varLambda$ additionally fulfils property~(Alg), then it only contains affinely regular $m$-gons for finitely many values of $m$.
\end{theorem}
\begin{proof}
For (i) $\Rightarrow$ (ii), let $P$ be an affinely regular $m$-gon in
$\varLambda$. There is then a non-singular affine transformation $\Psi
\!:\, \LR^2 \rightarrow \LR^2$ with $\Psi(R_{m}) = P$, where $R_{m}$ is
the regular $m$-gon with vertices given in complex form by $1, \zetam,
\dots, \zetam^{m - 1}$. If $m\in\{3,4,6\}$, condition (ii)
 holds trivially. Suppose $6\neq m\geq 5$. The pairs $\{1,\zetam\}$, $\{\zetam^{-1},\zetam^{2}\}$ lie on parallel lines and so do their images under $\Psi$. Therefore,
$$\frac{\vert \zetam^{2} - \zetam^{-1} \vert}{\vert \zetam - 1 \vert} = \frac{\vert \Psi(\zetam^{2}) - \Psi(\zetam^{-1}) \vert}{\vert \Psi(\zetam) - \Psi(1) \vert}\,.$$ Moreover, since $\Psi(\zetam^{2}) - \Psi(\zetam^{-1})$ and $\Psi(\zetam) - \Psi(1)$ are elements of $\varLambda -\varLambda$ and since $\vert z\vert^2=z\bar{z}$ for $z\in\LC$, we get the relation $$(1 + \zetam + \bar{\zeta}_{m})^2 = (1 + \zetam + \zetam^{-1})^2 = \frac{\vert \zetam^{2} - \zetam^{-1} \vert^2}{\vert \zetam - 1 \vert^2} = \frac{\vert \Psi(\zetam^{2}) - \Psi(\zetam^{-1}) \vert^2}{\vert \Psi(\zetam) - \Psi(1) \vert^2} \in \mathbbm{k}_{\varLambda}\,.$$ The pairs $\{\zetam^{-1},\zetam\}$, $\{\zetam^{-2},\zetam^{2}\}$ also lie on parallel lines. An argument similar to that above yields $$ (\zetam + \bar{\zeta}_{m})^2 = (\zetam + \zetam^{-1})^2 = \frac{\vert \zetam^{2} - \zetam^{-2} \vert^2}{\vert \zetam - \zetam^{-1} \vert^2} \in \mathbbm{k}_{\varLambda}\,.$$ By subtracting these equations, one gets the relation $$2(\zetam + \bar{\zeta}_{m}) + 1 \in \mathbbm{k}_{\varLambda}\,,$$ whence $\zetam + \bar{\zeta}_{m} \in \mathbbm{k}_{\varLambda}$, the latter being equivalent to the inclusion of the fields $ \mathbbm{k}_{m} \subset \mathbbm{k}_{\varLambda}$.

For (ii) $\Rightarrow$ (i), let $R_{m}$ again be the regular $m$-gon as defined in step (i) $\Rightarrow$ (ii). Since $m\geq 3$, the set $\{1,\zetam\}$ is an $\R$-basis of $\C$. Since $\varLambda$ is non-degenerate, there is an element $\tau\in\mathbbm{K}_{\varLambda}$ with non-zero imaginary part. Hence, one can define an $\R$-linear map $L\! :\, \LR^2\rightarrow \LR^2$ as the linear extension of $1 \mapsto 1$ and $\zetam \mapsto \tau$. Since $\{1,\tau\}$ is an $\R$-basis of $\C$ as well, this map is non-singular. Since $\mathbbm{k}_{m}\subset\mathbbm{k}_{\varLambda}$ and since $\{1,\zetam\}$ is a $\mathbbm{k}_{m}$-basis of $\mathbbm{K}_{m}$ (cf. Corollary~\ref{cr5}), the vertices of $L(R_{m})$, i.e., $L(1), L(\zetam), \dots, L(\zetam^{m - 1})$, lie in $\mathbbm{K}_{\varLambda}$, whence $L(R_{m})$ is a polygon in $\mathbbm{K}_{\varLambda}$. By property~(Aff), there is a non-singular affine transformation $\Psi\! :\, \LR^2 \rightarrow \LR^2$ such that $\Psi(L(R_{m}))$ is a polygon in $\varLambda$. Since compositions of non-singular affine transformations are non-singular affine transformations again, $\Psi(L(R_{m}))$ is an affinely regular $m$-gon in $\varLambda$.

For the additional statement, note that, since $\varLambda$ has property~(Alg), one has $[\mathbbm{k}_{\varLambda}:\LQ]<\infty$ by Remark~\ref{mrs}. Thus, $\mathbbm{k}_{\varLambda}/\LQ$ has only finitely many intermediate fields. The assertion now follows immediately from condition (ii) in conjunction with Corollary~\ref{unique}, Remark~\ref{unrem} and Lemma~\ref{unique2}.
\end{proof}

Let $\mathbbm{L}$ be an imaginary algebraic number field with $\overline{\mathbbm{L}}=\mathbbm{L}$ and let $\mathcal{O}_{\mathbbm{L}}$ be the ring of integers in $\mathbbm{L}$. Then, every translate $\varLambda$ of $\mathbbm{L}$ or $\mathcal{O}_{\mathbbm{L}}$ is non-degenerate and satisfies the properties (Alg) and (Aff). To this end, we first show that in both cases one has $\mathbbm{K}_{\varLambda}=\mathbbm{L}$. If $\varLambda$ is a translate of $\mathbbm{L}$, this follows immediately from the calculation
$$
\mathbbm{K}_{\varLambda}=\mathbbm{K}_{\mathbbm{L}}=\LQ(\mathbbm{L}\cup\overline{\mathbbm{L}})=\mathbbm{L}\,.
$$
If $\varLambda$ is a translate of $\mathcal{O}_{\mathbbm{L}}$, one has to observe that
$$
\mathbbm{K}_{\varLambda}=\mathbbm{K}_{\mathcal{O}_{\mathbbm{L}}}=\LQ(\mathcal{O}_{\mathbbm{L}}\cup\overline{\mathcal{O}_{\mathbbm{L}}})=\mathbbm{L}\,,
$$
 since
 $\overline{\mathbbm{L}}=\mathbbm{L}$ implies $\overline{\mathcal{O}_{\mathbbm{L}}}=\mathcal{O}_{\mathbbm{L}}$. In the first case, property (Aff) is evident, whereas, if $\varLambda$ is a translate of $\mathcal{O}_{\mathbbm{L}}$, property~(Aff) follows from the fact that there is always  a $\LZ$-basis of $\mathcal{O}_{\mathbbm{L}}$ that is simultaneously a $\LQ$-basis of $\mathbbm{L}$. Thus, if $F\subset \mathbbm{L}$ is a finite set, then a suitable translate of $aF$ is contained in $\varLambda$, where $a$ is defined as the least common multiple of the denominators of the $\LQ$-coordinates of the elements of $F$ with respect to a $\LQ$-basis of $\mathbbm{L}$ that is simultaneously a $\LZ$-basis of $\mathcal{O}_{\mathbbm{L}}$. Hence, for these two examples, property (Aff) may be replaced by the stronger property
\begin{itemize}
\item[(Hom)]
For all $F\in\mathcal{F}(\mathbbm{K}_{\varLambda})$, there is a
 homothety $h\!:\, \mathbbm{R}^{2} \rightarrow
\mathbbm{R}^{2}$ such that $h(F)\subset \varLambda$\,.
\end{itemize}

Thus, we have obtained the following consequence of Theorem~\ref{charac}.

\begin{coro}\label{numcor}
Let $\mathbbm{L}$ be an imaginary algebraic number field with $\overline{\mathbbm{L}}=\mathbbm{L}$ and let $\mathcal{O}_{\mathbbm{L}}$ be the ring of integers in $\mathbbm{L}$. Let $\varLambda$ be a translate of $\mathbbm{L}$ or a translate of $\mathcal{O}_{\mathbbm{L}}$. Further, let $m\in\mathbbm{N}$ with $m\geq 3$. Denoting the maximal real subfield of $\mathbbm{L}$ by $\mathbbm{l}$, the following statements are equivalent:
\begin{itemize}
\item[(i)]
There is an affinely regular $m$-gon in $\varLambda$.
\item[(ii)]
$\mathbbm{k}_{m}\subset\mathbbm{l}$.
\end{itemize}
Further, $\varLambda$ only contains affinely regular $m$-gons for finitely many values of $m$.
\end{coro}
 
\begin{rem}
In particular, Corollary~\ref{numcor} applies to translates of imaginary cyclotomic fields and their rings of integers, with $\mathbbm{l}=\kn$ for a suitable $n\geq 3$; cf. Facts~\ref{gau} and~\ref{p1} and also compare the equivalences of Corollary~\ref{th1mod} below. 
\end{rem}

\section{Application to cyclotomic model sets}

\begin{figure}
\centerline{\epsfysize=0.57\textwidth\epsfbox{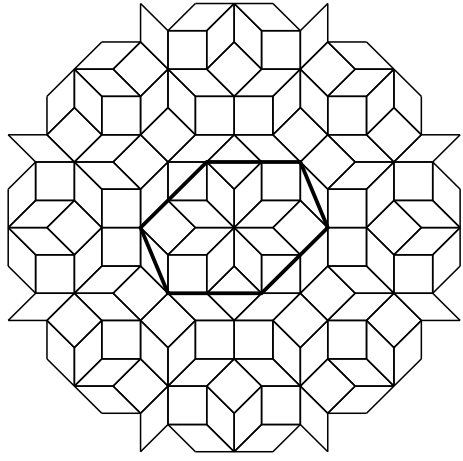}}
\caption{A central patch of the eightfold symmetric Ammann-Beenker tiling of the plane with squares and rhombi, both having edge length $1$. Therein, an affinely regular $6$-gon is marked.}
\label{fig:ab2}
\end{figure}

Remarkably, there are Delone subsets of the plane satisfying properties~(Alg) and (Hom). These sets were introduced as \emph{algebraic Delone sets} in~\cite[Definition 4.2]{H}. Note that algebraic Delone sets are always non-degenerate, since this is true for all relatively dense subsets of the plane. It was shown in~\cite[Proposition 4.15]{H} that the so-called \emph{cyclotomic model sets} $\varLambda$ are examples of algebraic Delone sets; cf. Section~4 of~\cite{H} and~\cite[Definition 4.6]{H} in particular for the definition of cyclotomic model sets. Any cyclotomic model set $\varLambda$ is contained in a translate of $\OO_n$, where $n\geq 3$, in which case the $\LZ$-module $\OO_n$ is called the \emph{underlying $\LZ$-module} of $\varLambda$. With the exception of the crystallographic cases of translates of the square lattice $\OO_4$ and translates of the triangular lattice $\OO_3$, cyclotomic model sets are aperiodic (i.e., they have no translational symmetries) and have long-range order; cf.~\cite[Remarks 4.9 and 4.10]{H}. Well-known examples of cyclotomic model sets with underlying $\LZ$-module $\OO_n$ are the vertex sets of aperiodic tilings of the plane like the Ammann-Beenker tiling~\cite{am,bj,ga} ($n=8$), the T\"ubingen triangle tiling~\cite{bk1,bk2} ($n=5$) and the shield tiling~\cite{ga} ($n=12$); cf.~\cite[Example 4.11]{H} for a definition of the vertex set of the Ammann-Beenker tiling and see~Figure~\ref{fig:ab2} and~\cite[Figure 1]{H} for illustrations. For further details and illustrations of the examples of cyclotomic model sets mentioned above, we refer the reader to~\cite[Section 1.2.3.2]{H2}. Clearly, any cyclotomic model set $\varLambda$  with underlying $\LZ$-module $\OO_n$ satisfies
\begin{equation}\label{eq1}
\mathbbm{K}_{\varLambda}\subset\LQ(\OO_n\cup\overline{\OO_n})= \Kn\,,\end{equation}
whence $\mathbbm{k}_{\varLambda}\subset\kn$. It was shown in ~\cite[Lemma 4.14]{H} that cyclotomic model sets $\varLambda$ with underlying $\LZ$-module $\OO_n$ even satisfy the following stronger version of property (Hom) above:

\begin{itemize}
\item[(\underline{Hom})]
For all $F\in\mathcal{F}(\mathbbm{K}_n)$, there is a
 homothety $h\!:\, \mathbbm{R}^{2} \rightarrow
\mathbbm{R}^{2}$ such that $h(F)\subset \varLambda$\,.
\end{itemize}

This property enables us to prove the following characterization. 

\begin{coro}\label{th1mod}
Let $m,n\in \mathbbm{N}$ with $m,n\geq 3$. Further, let $\varLambda$ be a cyclotomic model set with underlying $\LZ$-module $\OO_n$. The following statements are equivalent:
\begin{itemize}
\item[(i)]
There is an affinely regular $m$-gon in $\varLambda$.
\item[(ii)]
$\mathbbm{k}_{m}\subset\mathbbm{k}_{n}$.
\item[(iii)]
$m \in \{3,4,6\}$, or $\mathbbm{K}_{m}\subset \mathbbm{K}_{n}$.
\item[(iv)]
$m \in \{3,4,6\}$, or $m|n$, or $m=2d$ with $d$ an odd divisor of $n$.
\item[(v)]
$m \in \{3,4,6\}$, or $\OO_{m}\subset \OO_{n}$.
\item[(vi)]
$\oo_{m}\subset\oo_{n}$.
\end{itemize}
\end{coro}
\begin{proof}
Direction (i) $\Rightarrow$ (ii) is an immediate consequence of Theorem~\ref{charac} in conjunction with Relation~\eqref{eq1}. For direction (ii) $\Rightarrow$ (i), let $R_{m}$ again be the regular $m$-gon as defined in step (i) $\Rightarrow$ (ii) of Theorem~\ref{charac}. Since $m,n\geq 3$, the sets $\{1,\zetam\}$  and $\{1,\zetan\}$ are $\R$-bases of $\C$. Hence, one can define an $\R$-linear map $L\! :\, \LR^2\rightarrow \LR^2$ as the linear extension of $1 \mapsto 1$ and $\zetam \mapsto \zetan$. Clearly, this map is non-singular. Since $\mathbbm{k}_{m}\subset\mathbbm{k}_{n}$ and since $\{1,\zetam\}$ is a $\mathbbm{k}_{m}$-basis of $\mathbbm{K}_{m}$ (cf. Corollary~\ref{cr5}), the vertices of $L(R_{m})$, i.e., $L(1), L(\zetam), \dots, L(\zetam^{m - 1})$, lie in $\mathbbm{K}_{n}$, whence $L(R_{m})$ is a polygon in $\mathbbm{K}_{n}$. Because $\varLambda$ has property~(\underline{Hom}), there is a homothety $h\! :\, \LR^2 \rightarrow \LR^2$ such that $h(L(R_{m}))$ is a polygon in $\varLambda$. Since homotheties are non-singular affine transformations, $h(L(R_{m}))$ is an affinely regular $m$-gon in $\varLambda$. As an immediate consequence of Lemma~\ref{unique2}(b), we get the equivalence (ii) $\Leftrightarrow$ (iii). Conditions (iii) and (iv) are equivalent by Lemma~\ref{incl}. Finally, the equivalences (iii) $\Leftrightarrow$ (v) and (ii) $\Leftrightarrow$ (vi) follow immediately from Fact~\ref{p1}. 
\end{proof}

Although the equivalence (i) $\Leftrightarrow$ (iv) in Corollary~\ref{th1mod} is fully satisfactory, the following consequence deals with the two cases where condition (ii) can be used more effectively.

\begin{coro}\label{cor2}
Let $m,n\in \mathbbm{N}$ with $m,n\geq 3$. Further, let $\varLambda$ be a cyclotomic model set with underlying $\LZ$-module $\OO_n$. Consider $\phi$ on $\{n\in \mathbbm{N}\, |\, n\; \not\equiv \;2  \;(\operatorname{mod} 4)\}$. Then, one has:
\begin{itemize}
\item[(a)] If $n\in\{3,4\}$, there is an affinely regular $m$-gon in $\varLambda$ if and only if $m \in \{3,4,6\}$.  
\item[(b)] If $n\in\mathcal{S}$, there is an affinely regular $m$-gon in $\varLambda$ if and only if 
$$\left\{
\begin{array}{ll}
m \in \{3,4,6,n\}, & \mbox{if $n=8$ or $n=12$,}\\
m \in \{3,4,6,n,2n\}, & \mbox{otherwise.}
\end{array}\right.
$$ 
\end{itemize}
\end{coro}
\begin{proof}
By Lemma~\ref{phin2p}(a), $n\in\{3,4\}$ is equivalent to $\phi(n)/2=1$, thus condition (ii) of Corollary~\ref{th1mod}
specializes to $\mathbbm{k}_{m} = \Q$, the latter being equivalent to
$\phi(m)=2$, which means $m\in\{3,4,6\}$;  cf. Corollary \ref{cr5}. This
proves the part~(a). 

By Lemma~\ref{phin2p}(b), $n\in\mathcal{S}$ is equivalent to $\phi(n)/2\in\mathcal{P}$. By Corollary \ref{cr5}, this shows that $[\mathbbm{k}_{n}:\Q]=
\phi(n)/2 \in \mathcal{P}$. Hence, by condition (ii) of Corollary~\ref{th1mod}, one either gets
$\mathbbm{k}_{m}=\Q$ or $\mathbbm{k}_{m}=\mathbbm{k}_{n}$. The former
case implies $m\in\{3,4,6\}$ as in the proof of the part (a), while the proof
 follows from Lemma~\ref{unique2}(a) in conjunction with Corollary~\ref{unique} in the latter case.
\end{proof}

\begin{ex}
As mentioned above, the vertex set $\varLambda_{\rm AB}$ of the Ammann-Beenker tiling is a cyclotomic model set with underlying $\LZ$-module $\OO_8$. Since $8\in\mathcal{S}$, Corollary~\ref{cor2} now shows that there is an affinely regular $m$-gon in $\varLambda_{\rm AB}$ if and only if $m\in\{3,4,6,8\}$; see Figure~\ref{fig:ab2} for an affinely regular $6$-gon in $\varLambda_{\rm AB}$. The other solutions are rather obvious, in particular the patch shown also contains the regular $8$-gon $R_8$, given by the $8$th roots of unity. 
\end{ex}

For further illustrations and explanations of the above results, we refer the reader to~\cite[Section 2.3.4.1]{H2} or~\cite[Section 5]{H0}. This references also provide a detailed description of the construction of affinely regular $m$-gons in cyclotomic model sets, given that they exist.

\section{Application to discrete tomography of cyclotomic model sets}

\emph{Discrete tomography} is concerned with the  
inverse problem of retrieving information about some finite 
object from information about its slices; cf.~\cite{G,GG,H0,H,H2} and also see the refences therein. A typical example is the \emph{reconstruction} of a finite point set
from its ({\em
  discrete parallel}\/) {\em $X$-rays} in a small number of directions. In the following, we restrict ourselves to the planar case. 

\begin{defi}\label{xray..}
Let $F\in
\mathcal{F}(\mathbbm{R}^{2})$, let $u\in
\mathbb{S}^{1}$ be a direction and let $\mathcal{L}_{u}$ be the set
of lines in direction $u$ in $\mathbbm{R}^{2}$. Then, the ({\em
  discrete parallel}\/) {\em X-ray} of $F$ {\em in direction} $u$ is
the function $X_{u}F: \mathcal{L}_{u} \rightarrow
\mathbbm{N}_{0}:=\mathbbm{N} \cup\{0\}$, defined by $$X_{u}F(\ell) :=
\operatorname{card}(F \cap \ell\,) =\sum_{x\in \ell}
\mathbbm{1}_{F}(x)\,.$$ 
\end{defi}

In~\cite{H}, we studied the problem of \emph{determining} convex subsets of algebraic Delone sets $\varLambda$ by $X$-rays. Solving this problem amounts to find small sets $U$ of suitably prescribed $\varLambda$-directions with the property that different convex subsets of $\varLambda$ cannot have the same $X$-rays in the directions of $U$. More generally, one defines as follows.

\begin{defi}
Let $\mathcal{E}\subset
\mathcal{F}(\mathbbm{R}^{2})$, and let $m\in\LN$. Further, let $U\subset\mathbb{S}^{1}$ be a
finite set of directions. We say that $\mathcal{E}$ is {\em determined} by the $X$-rays in the directions of $U$ if, for all $F,F' \in \mathcal{E}$, one has
$$
(X_{u}F=X_{u}F'\;\,\forall u \in U) \;  \Longrightarrow\; F=F'\,.
$$
\end{defi}

Let $\varLambda\subset\R^2$ be a Delone set and let $U\subset \mathbb{S}^{1}$ be a set of two or more pairwise non-parallel $\varLambda$-directions. Suppose the existence of a $U$-polygon $P$
in $\varLambda$. Partition the vertices of $P$ into two
disjoint sets $V,V'$, where the elements of these sets alternate round
the boundary of $P$. Since $P$ is a $U$-polygon, each line
in the plane parallel to some $u\in U$ that contains a point in $V$
also contains a point in $V'$. In particular, one sees that $\operatorname{card}(V)=\operatorname{card}(V')$. Set
$$
C:=(\varLambda\cap P)\setminus (V\cup V')
$$
and, further, $F:=C\cup V$ and $F':=C\cup V'$. Then, $F$ and $F'$ are different convex subsets of $\varLambda$ with the same $X$-rays in the directions of $U$. We have just proven direction (i) $\Rightarrow$ (ii) of the following equivalence, which particularly applies to cyclotomic model sets, since any cyclotomic model set is an algebraic Delone set by~\cite[Proposition 4.15]{H}.

\begin{theorem}\cite[Theorem 6.3]{H}\label{characungen}
Let $\varLambda$ be an algebraic Delone set and let $U\subset \mathbb{S}^{1}$ be a set of two or more pairwise non-parallel $\varLambda$-directions. The following statements are equivalent:
\begin{itemize}
\item[(i)]
$\mathcal{C}(\varLambda)$ is determined by the $X$-rays in the directions of $U$.
\item[(ii)]
There is no $U$-polygon in $\varLambda$.
\end{itemize}
\end{theorem}

\begin{rem}\label{trivrem}
Trivially, any affinely regular $m$-gon $P$ in $\varLambda$ with $m$ even  is a $U$-polygon in $\varLambda$ with respect to any set $U\subset\mathbb{S}^1$ of pairwise non-parallel directions having the property that each element of $U$ is parallel to one of the edges of $P$. The set $U$ then consists only of $\varLambda$-directions and, moreover, satisfies $\operatorname{card}(U)\leq m/2$. 
\end{rem}

By combining Corollary~\ref{th1mod}, direction (i) $\Rightarrow$ (ii) of Theorem~\ref{characungen} and Remark~\ref{trivrem}, one immediately obtains the following consequence.

\begin{coro}\label{corofin}
Let $n\geq 3$ and let $\varLambda$ be a cyclotomic model set with underlying $\LZ$-module $\OO_n$. Suppose that there exists a natural number $k\in\LN$ such that, for any set $U$ of $k$ pairwise non-parallel $\varLambda$-directions, the set $\mathcal{C}(\varLambda)$ is determined by the $X$-rays in the directions of $U$. Then, one has
$$
k>\,\max\left\{3,\tfrac{\operatorname{lcm}(n,2)}{2}\right\}\,.
$$
\end{coro}

\begin{rem}
In the situation of Corollary~\ref{corofin}, the question of existence of a suitable number $k\in\LN$ is a much more intricate problem. So far, it has only been answered affirmatively by Gardner and Gritzmann in the case of translates of the square lattice ($n=4$), whence corresponding results hold for all translates of planar lattices, in particular for translates of the triangular lattice ($n=3$); cf.~\cite[Theorem 5.7(ii) and (iii)]{GG}. More precisely, it is shown there that, for these cases, the number $k=7$ is the smallest among all possible values of $k$. It would be interesting to know if suitable numbers $k\in\LN$ exist for all cyclotomic model sets. 
\end{rem}

Let us finally note the following relation between $U$-polygons and affinely regular polygons. The proof uses a beautiful theorem of Darboux~\cite{D} on second midpoint polygons; cf.~\cite{GM} or~\cite[Ch. 1]{G}.

\begin{prop}\cite[Proposition 4.2]{GG}\label{uaffine}
If $U\subset \mathbb{S}^{1}$ is a finite set of directions, there exists a $U$-polygon if and only if there is an affinely regular polygon such that each direction of $\,U$ is parallel to one of its edges.
\end{prop}

\begin{rem}
A $U$-polygon need not itself be an affinely regular polygon, even if it is a $U$-polygon in a cyclotomic model set; cf.~\cite[Example 4.3]{GG} for the case of planar lattices and~\cite[Example 2.46]{H2} or~\cite[Example 1]{H0} for related examples in the case of aperiodic cyclotomic model sets. 
\end{rem}

\section*{Acknowledgements}
I am indebted to Michael Baake, Richard J. Gardner, Uwe Grimm and Peter A. B. Pleasants for their cooperation and for useful hints on the manuscript. Interesting discussions with Peter Gritzmann and Barbara Langfeld are gratefully acknowledged.

\end{document}